%% file: multipliers-arxiv.tex
\input neumac
\input multipliers-arxiv.num
\nicknum=0
\def\tr{{\rm tr}}\def\sgn{{\rm sgn\thinspace}}
\nopagenumbers
\let\header=N
\immediate\newwrite\num\immediate\openout\num=multipliers-arxiv.num
\def\leftheadline{\ninepoint\folio\hfill
Multiplier systems}%
\def\transpose#1{\kern1pt{^t\kern-1pt#1}}%

\noindent
{\titelfont Multiplier systems for Siegel modular groups}%
\vskip 1.5cm
\centerline{Eberhard Freitag, mathematisches Institut, Universit"at Heidelberg}
\centerline{Adrian Hauffe-Waschb"usch, Institut f"ur Mathematik, RWTH Aachen}
\centerline{2020}
\vskip1cm \noindent
\noindent
\centerline{\vbox{\noindent\hsize=10cm
{\ninepoint{\bf Abstract}
\smallni
Deligne proved in [De] (s.\ also [Hi], 7.1) that the weights of
Siegel modular forms on any congruence subgroup
of the Siegel modular group of genus $g>1$ must be integral or half integral. Actually he proved that
for a system $v(M)$ of complex numbers of absolute value 1
$$v(M)\det(CZ+D)^r\qquad(r\in\rz)$$
can be an automorphy factor only if $2r$ is integral. We give a different proof for this.
It uses Mennicke's result that subgroups of finite index of the Siegel modular group
are congruence subgroups and some techniques from the paper [BMS] of Bass-Milnor-Serre.
}}}
\vskip 1cm\noindent
{\paragratit Introduction}%
\medskip\noindent
We fix a natural number $g$ (which later will be 2). We denote by $E=E^{(g)}$ the $g\times g$-unit matrix
and by
$$I=I^{(g)}=\pmatrix{0&-E\cr E&0}$$
the standard alternating matrix.
The symplectic group $\Sp(g,\rz)$ consists of all
$M\in\GL(2g,\rz)$ with the property $M'IM=I$. Here $M'$ denotes the transposed matrix of $M$.
We consider the usual action $MZ=(AZ+B)(CZ+D)^{-1}$ of the real
symplectic group $\Sp(g,\rz)$ on the Siegel upper half plane.
The function
$$J(M,Z)=\det(CZ+D)$$
has no zeros on the half plane.
Since the half plane is convex, there exists a continuous choice
$L(M,Z)=\arg J(M,Z)$ of the argument. We normalize it such that
it is the principal value for $Z=\imag E$  where $E$ denotes the
unit matrix.
Recall that the principal value $\Arg(a)$ is defined such that it
is in the interval $(-\pi,\pi]$.
So we have
$$L(M,\imag E)=\Arg(J(M,\imag E))\in(-\pi,\pi].$$
We consider
$$w(M,N):={1\over 2\pi}\bigl((L(MN,Z)-L(M,NZ)-L(N,Z)\bigr).$$
Obviously,
$$e^{2\pii w(M,N)}=1.$$ Hence $w(M,N)$ is a constant (independent of $Z$),
$$w(M,N)\in \gz.$$
\Proclaim
{Remark}
{The function $w:\Sp(n,\rz)\times\Sp(n,\rz)\to\gz$ is a cocycle
in the following sense:
$$\eqalign{w(M_1M_2,M_3)+w(M_1,M_2)&=w(M_1,M_2M_3)+w(M_2,M_3),\cr
w(E,M)=w(M,E)&=0.\cr}$$
}
\finishproclaim
The computation of $w(M,N)$ in genus 1 is easy for the following reason.
From the definition we have
$$2\pi w(M,N)=\Arg((c\alpha+d\gamma)\imag+c\beta+d\gamma)-\arg(c N(\imag)+d)-\Arg(\gamma\imag+\delta)$$
for
$$M=\pmatrix{a&b\cr c&d},\quad N=\pmatrix{\alpha&\beta\cr\gamma&\delta}$$
where $\arg(cN(\imag)+d)$ is obtained from the principal value of $\arg (c\imag +d)$ through
continuous continuation. But $cz+d$ for $z$ in the upper half plane
never crosses the real axis. Hence the
result of the continuation is the principal value too. So all three arguments in the definition
of $w(M,N)$ are the principal values (in genus 1). This makes it easy to compute $w$.
We rely on tables for the values of $w$ which have been derived by Petersson and reproduced by Maass [Ma1],
Theorem 16.
\proclaim
{Lemma}
{Let $M={*\;\;\;*\choose m_1\,m_2}$, $S={a\,b\choose c\,d}$ be two real matrices with determinant $1$ and
$(m_1',m_2')$ the second row of the matrix $MS$. Then
$$4w(M,S)=\cases{\sgn c+\sgn m_1-\sgn m_1'-{{\rm sgn}}(m_1cm_1')& if $m_1cm_1'\ne 0$,\cr
-(1-\sgn c)(1-\sgn m_1)& if $cm_1\ne0, m_1'=0$,\cr
(1+\sgn c)(1-\sgn m_2)&if $cm_1'\ne 0,m_1=0$,\cr
(1-\sgn a)(1+\sgn m_1)&if $m_1m_1'\ne 0,c=0$,\cr
(1-\sgn a)(1-\sgn m_2)&if $c=m_1=m_1'=0$.\cr}$$
{\bf Corollary.}
Assume that $m_1cm_1'\ne 0$ and that $m_1m_1'>0$ or $m_1c<0$. Then $w(M,S)=0$.}
MP%
\finishproclaim
We give an example.
\proclaim
{Lemma}
{We have
$$w\biggl(\pmatrix{a&b\cr c&d},\pmatrix{1&x\cr0&1}\biggr)=
w\biggl(\pmatrix{1&x\cr0&1},\pmatrix{a&b\cr c&d}\biggr)=0.$$
}
TraTr%
\finishproclaim
We denote by $\Gamma_g[q]$ the principal congruence subgroup level $q$.
This is the kernel of the natural homomorphism
$\Sp(g,\gz)\to\Sp(g,\gz/q\gz)$.
\neupara{Some special values of the cocycle}%
We give some examples for values of $w$ in genus $g>1$.
\proclaim
{Lemma}
{One has
$$w\biggl(\pmatrix{E&S\cr0&E},M\biggr)=0.$$}
LTra%
\finishproclaim
The proof is trivial and can be omitted.\qed
\smallskip
\proclaim
{Lemma}
{Let $g=2$ and
$$P=\pmatrix{0&1&0&0\cr1&0&0&0\cr 0&0&0&1\cr0&0&1&0}.$$
We have
$$w(P,M)=w(M,P)=\cases{0&if $\,\im\det(\imag C+D)<0$,\cr -1& if $\,\im\det(\imag C+D)>0$.\cr}$$}
Pval%
\finishproclaim
{\it Proof.} Let $z:=\det(C\imag+D)$. One computes
$$2\pi w(P,M)=2\pi w(M,P)=\Arg (-z)-\Arg(z)-\Arg(-1).\eqno\square$$
\proclaim
{Definition}
{The {\emph Siegel parabolic group} consists of all symplectic matrices of the form
$$\pmatrix{A&B\cr0&D}.$$
The two {\emph Klingen parabolic groups} in the case $g=2$ consist of all symplectic matrices of the form
$$\pmatrix{a_1&0&b_1&b_2\cr a_3&a_4&b_3&b_4\cr
c_1&0&d_1&d_2\cr0&0&0&d_4}\quad\hbox{resp.}\quad
\pmatrix{a_1&a_2&b_1&b_2\cr 0&a_4&b_3&b_4\cr
0&0&d_1&0\cr0&c_4&d_3&d_4}.$$
}
DefSK%
\finishproclaim
There is a character on the Siegel parabolic group
$$\varepsilon\pmatrix{A&B\cr 0&D}=\det(D).$$
For an element $M$ of the Siegel parabolic group, the expression
$\det(CZ+D)=\det(D)$ is independent of $Z$. Hence
$$L(M,Z)=0\quad\hbox{if}\ \varepsilon(M)>0.$$
An immediate consequence is the following lemma.
\proclaim
{Lemma}
{For two elements $M,N$ of the Siegel parabolic group we have
$w(M,N)=0$ if $\varepsilon(M)>0$.}
ParM%
\finishproclaim
\proclaim
{Lemma}
{Let $g=2$ and let $M$ be a Klingen parabolic matrix and $N$ a Siegel parabolic matrix with $\varepsilon(N)>0$. Then
$w(M,N)=0$.}
KSz%
\finishproclaim
{\it Proof.} Since $L(N,\imag E)=1$, we have to show that the arguments of of
$J(MN,\imag E)$ and of $L(M,N(\imag E))$ are the same.  Both determinants are equal. But the argument of
the first is the principal part and that of the second is defined by continuation from the argument of
$J(M,\imag E)$. Hence it is sufficient to show that the principal part of the argument of
$L(M,Z)$ is continuous. This is the case if $\im J(M,Z)$ is always $\ge 0$ or always $< 0$. Actually,
for the first Klingen parabolic group
$$\im J(M,Z)=c_1d_4\im z_0\qquad\hbox{where}\quad Z=\pmatrix{z_0&*\cr*&*}.$$
The argument for the second Klingen parabolic group is the same.
This proves the lemma.\qed
\proclaim
{Lemma}
{Let $g=2$ and let
$$M=\pmatrix{E&S\cr0&E}.$$
Then
$$w(I,M)=\cases{0&if $\tr(S)\ge 0$,\cr -1& else.\cr}$$}
ITra%
\finishproclaim
{\it Proof.}
From the definition we have
$$2\pi w(I,M)=\Arg\det(\imag E+S)-\Arg\det(E)-\arg\det(\imag E+S).$$
The third argument is defined through  continuation
along $\det(\imag E+tS)$, beginning from $t=0$ to $t=1$. For $t=0$ we have to
take the principal value which is $\pi$. The imaginary part of $\det(\imag E+tS)$
equals $t\tr(S)$. In the case $\tr(S)\ge 0$ we keep the principal value.
But if it is negative we make a jump by $-2\pi$.\qed
\proclaim
{Lemma}
{Let $g=2$ and let
$$M=\pmatrix{E&0\cr S&E}.$$
Then
$$w(M,I)=\cases{-1&if $\tr(S)\ge 0$,\cr 0& else.\cr}$$}
TraI%
\finishproclaim
{\it Proof.} Let $z=\det(\imag S+E)$. One computes $2\pi w(M,I)=\Arg(-z)-\pi-\Arg(z)$.
This depends on the imaginary part of $z$ which is $\tr(S)$.\qed
\neupara{Multipliers}%
\proclaim
{Definition}
{Let $\Gamma\subset\Sp(g,\rz)$ be an arbitrary subgroup
and let $r$ be a real number. A system $v(M)$, $M\in\Gamma$, of complex numbers of absolute value $1$
is called a multiplier system of weight $r$ if
$$v(MN)\equiv v(M)v(N)\sigma(M,N)$$
where
$$\sigma(M,N)=\sigma_r(M,N):=e^{2\pii r w(M,N)}.$$
}
TriVw%
\finishproclaim
The elliptic modular group
$\Sp(1,\gz)=\SL(2,\gz)$ admits  multipliers for every real $r$.
One can construct them by means of the discriminant function $\Delta$.
This is a modular form without zeros. Hence we can choose a
holomorphic power $f(z)=\Delta(z)^{r/12}$. This can be used to
construct a multiplier.
\smallskip
Maass [Ma2] proved that the full Siegel modular group of genus $g>1$ admits only multipliers for
integral $r$ and their values can be only $\pm1$.
As a consequence (s.~[Ch]), for every multiplier system on a subgroup $\Gamma$ of finite index of the modular group
the weight $r$ is rational and the values of $v$ are contained in a finite subgroup of $S^1$.
\smallskip
Let $\Gamma_{g,\vartheta}$ be the theta group of degree $g$.
It consists of all integral symplectic matrices such that
$AB'$ and $CD'$ have even diagonal. The function
$$\vartheta(Z)=\sum_{n\in\gz^g}e^{2\pii n'Zn}$$
is a modular form of weight $1/2$ on the theta group. It can
be used to construct a multiplier system of weight $1/2$.
\smallskip
The result of Deligne states:
\proclaim
{Theorem}
{Let $g>1$ and let $\Gamma\subset\Sp(g,\gz)$ be any
subgroup of finite index of the Siegel modular group.
Multiplier systems of weight $r$ can only exist if $2r$ is integral.}
ThD%
\finishproclaim
It is sufficient to prove this in the case $g=2$. So we assume from now on $g=2$.
\smallskip
We assume that a natural number $q'$ is given and that
$v$ is a multiplier system of weight $r$ on
$\Gamma_2[q']$.
\smallskip
For any $L\in\Sp(2,\gz)$ we can consider a conjugate
multiplier system [FB] that is defined by
$$\tilde v(M)=v(LML^{-1}){\sigma(LML^{-1},L)\over \sigma(L,M)}.$$
It is easy to check that this is a multiplier system.
The quotient of two multiplier systems of the same weight is a homomorphism, as we know into a finite group.
Since every subgroup of finite index of the Siegel modular group is a congruence subgroup,
we obtain $\tilde v(M)=v(M)$ on some subgroup $\Gamma_2[q]\subset\Gamma_2[q']$ (where $q$ may depend on $L$).
\proclaim
{Lemma}
{For given $L$ in the full modular group there exists a multiple $q$ of $q'$ such that
such that
$$v(M)=v(LML^{-1}){\sigma(LML^{-1},L)\over\sigma(L,M)}$$
for each $M\in\Gamma_2[q]$.
}
EleO%
\finishproclaim
This will be used for several matrices, in particular for $M=I$. 
\proclaim
{Proposition}
{There exists a multiple $q$ of $q'$ such that the following holds.
Let $U$ be an element from the subgroup that is generated by the matrices
${1\,q\choose 0\,1}$ and ${1\,0\choose q\,1}$. Let
$M$ be a matrix from $\Gamma_2[q]$ of the form
$$M=\pmatrix{U'&*\cr0&U^{-1}}\quad\hbox{or}\quad M=\pmatrix{E&0\cr*&E}.$$
Then $v(M)=1$.}
ProO%
\finishproclaim
{\it Proof.} The matrices of the first type build a finitely generated group.
The number of generators is independent on $q$.
It is enough to prove $v(M)=1$ for the generators, since $w(M,N)=0$ for all
$M,N$ in this group. We also have $v(M)^n=v(M^n)$. Since the values of $v$ are
contained in a finite group, we find an $n$ such that $v(M^n)=1$ for all of
the generators.
\smallskip
The second case is more difficult. Due to Lemma \EleO\ it is sufficient
to prove
$\sigma(IMI^{-1},I)=\sigma(I,M)$ for translation matrices $M$. This follows from
the Lemmas \ITra\ and \TraI.\qed
\neupara{Embedded subgroups}%
We have to consider three embeddings of $\SL(2,\gz)$ into $\Sp(2,\gz)$, namely
$$\iota_1,\iota_2,\iota_3:\SL(2)\lo\Sp(2),$$
$$\iota_1\pmatrix{a&b\cr c&d}=
\pmatrix{a&0&b&0\cr 0&1&0&0\cr c&0&d&0\cr 0&0&0&1},\quad
\iota_2\pmatrix{a&b\cr c&d}=
\pmatrix{1&0&0&0\cr0&a&0&b\cr0&0&1&0\cr0&c&0&d},$$
$$\iota_3(M)=\pmatrix{M&0\cr0&M'^{-1}}=\pmatrix{a&b&0&0\cr c&d&0&0\cr 0&0&d&-c\cr0&0&-b&a}.$$
We have $w(\iota_3(M),\iota_3(N))=0$. Hence $M\mapsto v(\iota_3(M))$ is a homomorphism into a finite group.
Its kernel is a subgroup of finite index in $\SL(2,\gz)$. We will show that it is in fact a
congruence subgroup.
\smallskip
Let $P$ as in Lemma \Pval. We have
$$P\iota_1(M)P^{-1}=\iota_2(M).$$
From Lemma \Pval\ follows
$w(\iota_2(M),P)=w(P,\iota_1(M))$. Hence we obtain from Lemma \EleO\
the following result.
\proclaim
{Lemma}
{We have
$$v(\iota_1( M))=v(\iota_2(M))$$
for $M\in\Gamma_1[q]$.}
iEiZ%
\finishproclaim
For sake of simplicity we write
$$v(M)=v(\iota_1(M))=v(\iota_2(M)).$$
This is a multiplier system  in genus 1.
We have
$$w(M,N)=w(\iota_\nu(M),\iota_\nu(N)),\quad\hbox{for}\quad\nu=1,2.$$
\smallskip
\proclaim
{Lemma}
{The value $v(M)$, $M\in\Gamma_1[q]$, depends only on the second row of~$M$.}
SRd%
\finishproclaim
{\it Proof.} When $M,N$ have the same second row, then
${1\,x\choose0\,1}M=N$. We know $w\bigl({1\,x\choose0\,1},M\bigr)=0$
and $v{1\,x\choose0\,1}=1$ (Proposition \ProO).\qed
\proclaim
{Lemma}
{Assume that $v$ is a multiplier system of weight $r$ on $\Gamma_2[q']$. There exists a
multiple $q$ of $q'$ such that
for $M={a\,b\choose c\,d}\in\Gamma_1[q]$ we have
$$v\pmatrix{d_1&-c_1&0&0\cr-b_1&a&0&0\cr0&0&a&b_1\cr0&0&c_1&d_1}\cdot
v\pmatrix{a&0&b_2&0\cr0&1&0&0\cr c_2&0&d_2&0\cr0&0&0&1}
= v\pmatrix{1&0&0&0\cr0&a&0&b_1^2b_2\cr0&0&1&0\cr0&c_1^2c_2&0&y}$$
where
$$y=d_1-b_1c_1d_2+c_1c_2b_1b_2d_1.$$
}
ZweiV%
\finishproclaim
{\it Proof.}
The proof depends on a certain relation which occurs in [BMS] during the proof
of Lemma 13.3. We reproduce it here.
We set
$$\eqalign{&
H_1=\pmatrix{d_1&-c_1&0&0\cr-b_1&a&0&0\cr0&0&a&b_1\cr0&0&c_1&d_1},\
H_2=\pmatrix{a&0&b_2&0\cr0&1&0&0\cr c_2&0&d_2&0\cr0&0&0&1},\cr&
H_3=\pmatrix{1&0&0&0\cr0&a&0&b_1^2b_2\cr0&0&1&0\cr0&c_1^2c_2&0&y}.\cr}$$
We consider the matrices
$$\eqalign{
&R_1=\pmatrix{1&0&0&0\cr b_1&1&0&0\cr0&0&1&-b_1\cr0&0&0&1},\quad
R_2=\pmatrix{1&0&0&0\cr0&1&0&0\cr ac_2&c_1c_2&1&0\cr c_1c_2&0&0&1}\cr&
R_3=\pmatrix{1&c_1&0&0\cr0&1&0&0\cr0&0&1&0\cr0&0&-c_1&1},\quad
R_4=\pmatrix{1&0&-ad_1^2b_2&b_1b_2d_1\cr0&1&b_1b_2d_1&0\cr0&0&1&0\cr0&0&0&1\cr}
.\cr}$$
Now a direct computation shows
$$\eqalign{
R_2H_3=H_1H_2R_1R_3R_4.\cr}$$
We have to compute $w$-values. We assume that $c_1c_2\ne 0$. First we treat $w(R_2,H_3)$.
We have
$$R_2H_3=\pmatrix{*&*&*&*\cr*&*&*&*\cr ac_2&ac_1c_2&1&b_1^2b_2c_1c_2\cr
c_1c_2&c_1^2c_2&0&y}.$$
We are going to compute $w(R_2,H_3)$.
A direct computation gives
$$\im J(R_2H_3,\imag E)=c_2(1+c_1^2).$$
Next we treat $J(R_2,H_3(\imag E))$. Here  the argument has to be defined by continuation
from the principal value of the argument of $J(R_2,\imag E)$. We can do this along the
straight line from $\imag E$ to $H_3(\imag E))$. The points on this line are of the form
${\imag\,0\choose0\,\tau}$ where $\tau$ is in the upper half plane. One computes
$$J\Bigl(R_2,{\imag\,0\choose0\,\tau}\Bigr)=\det\pmatrix{1+ac_2\imag&c_1c_2\tau\cr c_1c_2\imag&1}.$$
The real part is $1+c_1^2c_2^2\im\tau$ which is positive. Hence the principal value of the argument
is continuous along the line. So we see
$$L(R_2,H_3(\imag E))\in (-\pi,\pi].$$
Finally we compute
$$\im J(H_3,\imag E)=c_1^2c_2.$$
Now we see that the imaginary part of $J(R_2H_3,\imag E)$ and $J(H_3,\imag E)$ have the same
sign (namely the sign of $c_2$). Hence their arguments are both contained
in $(0,\pi)$ or in $(-\pi,0)$. This means that
$2\pi w(R_2,H_3)$ is contained in $(0,\pi)-(-\pi,\pi]-(0,\pi)$ or in
$(-\pi,0)-(-\pi,\pi]-(-\pi,0)$. This is $(-2\pi,2\pi)$ in both cases. We obtain
$$w(R_2,H_3)=0.$$
The case $c_1c_2=0$ is easy and can be omitted.
\smallskip
From Lemma \KSz\ we can take $w(H_2,R_1R_3R_4)=0$. For trivial reason one has
$w(H_1,H_2R_1R_3R_4)=0$.
Now we evaluate
$$v(R_2H_3)=v(H_1H_2R_1R_3R_4).$$
The left hand side is
$$v(R_2)v(H_3)\sigma(R_2,H_3)=v(R_2)v(H_3).$$
But $v(R_2)=1$ (Proposition \ProO). Hence the left hand side is just $v(H_3)$.
The right hand side is
$$v(H_1)v(H_2R_1R_3R_4)\sigma(H_1,H_2R_1R_3R_4)=v(H_1)v(H_2R_1R_3R_4).$$
Similarly we see
$$v(H_2R_1R_2R_3)=v(H_2)v(R_1R_2R_3)\sigma(H_2,R_1R_3R_4)=v(H_2)v(R_1R_3R_4).$$
From Proposition \ProO\ we know $v(R_1R_3R_4)=1$.
Hence we get $v(H_3)=v(H_1)v(H_2)$.\qed
\neupara{Mennicke symbol}%
We have seen that $v(\iota_1(M))=v(\iota_2(M))$ depends only on the second row of $M\in\Gamma_1[q]$. Hence we
can define
$$\Bigl\{{c\atop d}\Bigr\}=v\pmatrix{a&0&b&0\cr 0&1&0&0\cr c&0&d&0\cr 0&0&0&1}^{-1}=
v\pmatrix{1&0&0&0\cr0&a&0&b\cr0&0&1&0\cr0&c&0&d}^{-1}.$$
We also can define
$$\Bigl[{b\atop a}\Bigr]=v\pmatrix{a&b&0&0\cr c&d&0&0\cr 0&0&d&-c\cr0&0&-b&a}.$$
It is clear that this does not depend on the choice of $c,d$.
\proclaim
{Proposition}
{For a suitable multiple $q>2$ of $q'$
the bracket $\bigl[{b\atop a}]$ is a Mennicke symbol. This means that it is a function
on the set of all coprime $(a,b)$ with the property $a\equiv 1\mod q$ and $b\equiv0\mod q$
such that the following properties hold.
\smallni
{\rm MS1} It is invariant under the transformations $(a,b)\mapsto (a+xb,b)$ and $(a,b)\mapsto (a,b+qay)$
for integral $x,y$.
\smallni
{\rm MS2} It satisfies the rule
$$\Bigl[{b_1b_2\atop a}\Bigl]=\Bigl[{b_1\atop a}\Bigl]\Bigl[{b_2\atop a}\Bigl].$$}
MenS%
\finishproclaim
{\it Proof of\/ {\rm MS1}.}
We notice that $w$ is trivial on the image of $\iota_3$. Hence $v$ is a character on this group.
The invariance under $(a,b)\mapsto (a,b+qay)$ follows from the equation
$$\pmatrix{a&b\cr c&d}\pmatrix{1&qy\cr 0&1}=\pmatrix{a&b+qay\cr*&*}.$$
To prove the invariance under $(a,b)\mapsto (a+xb,b)$, we consider
$$\pmatrix{1&0\cr-x&1}\pmatrix{a&b\cr c&d}\pmatrix{1&0\cr x&1}=\pmatrix{a+xb&b\cr*&*}.$$
Due to Lemma \EleO\ we can assume that $v(\iota_3(M))$ is invariant under conjugation
with $\iota_3{1\,0\choose 1\,1}$. This proves MS1.
\smallni
{\it Proof of\/ {\rm MS2}.}
The proof of MS2 needs two Lemmas which we now have to formulate and prove now.
We make use of
$$v(\iota_\nu(M^{-1}))=v(\iota_\nu(M))^{-1}.$$
This is true since in genus 1 one has $w(M,M^{-1})=0$. (This is a general rule for
$c\ne 0$ but also for $c=0$ and $a>0$. But in our case $c=0$ implies $a=1$ since
we assume $q>2$.)
This relation implies
$$\Bigl\{{c\atop d}\Bigl\}=\Bigl\{{-c\atop a}\Bigl\}^{-1}.$$
From Lemma \ZweiV\ we get after the replacement,
$c_2\mapsto-c_2$ the following general rule (compare Lemma 13.3 in [BMS]).
\proclaim
{Lemma}
{Let $a-1\equiv c_1\equiv c_2\equiv 0$ mod $q$ and let $a,c_1$ and $a,c_2$ be coprime. Then
$$\Bigl[{c_1\atop a}\Bigr]\Bigl\{{c_2\atop a}\Bigr\}=\Bigl\{{c_1^2c_2\atop a}\Bigr\}.$$
}
emaA%
\finishproclaim
We need also the following simple lemma.
\proclaim
{Lemma}
{We have
$$\Bigl\{{1-a\atop a}\Bigr\}=1$$
for $a\equiv 1\mod q$.}
emaB%
\finishproclaim
{\it Proof.} We use
$$\pmatrix{1&1\cr0&1}\pmatrix{1&0\cr a-1&1}\pmatrix{1&-1\cr0&1}=\pmatrix{2-a&a-1\cr 1-a&a}.\eqno\square$$
We insert un Lemma \emaA\ now $c_2=1-a$ to obtain the following formula.
$$\Bigl[{c\atop a}\Bigr]=\Bigl\{{c^2(1-a)\atop a}\Bigr\}.$$
Before we continue, we mention that $\{\}$ is not  a Mennicke symbol. It does not satisfy MS1.
\proclaim
{Lemma}
{We have
$$\Bigl\{{c\atop d}\Bigr\}=\Bigl\{{c\atop d+yc}\Bigr\}$$
and
$$\Bigl\{{c+xqd\atop d}\Bigr\}=\Bigl\{{c\atop d}\Bigr\}e^{2\pii rs}.$$
where
$$s=w\biggl(\pmatrix{*&*\cr c&d},\pmatrix{1&0\cr qx&1}\biggr).$$
}
GnM%
\finishproclaim
{\it Proof.}
The first relation
can be derived from
$$\pmatrix{1&-y\cr 0&1}\pmatrix{*&*\cr c&d}\pmatrix{1&y\cr0&1}=\pmatrix{*&*\cr c&d+cy}.$$
To derive the second one we consider the relation
$$\pmatrix{*&*\cr c&d}\pmatrix{1&0\cr qx&1}=\pmatrix{*&*\cr c+dxq&d}.$$
It shows
$$\Bigl\{{c+dxq\atop d}\Bigr\}=\Big\{{c\atop d}\Bigr\}\;e^{2\pii rs}.$$
The $w$-value $s$ is usually not zero.\qed
\smallni
{\bf Proof of Proposition \MenS\ (MS2) continued.}
Now we use
$$\pmatrix{*&*\cr c^2&a}\pmatrix{1&0\cr-c^2&1}=\pmatrix{*&*\cr c^2-ac^2&a}.$$
From the corollary of the table of Maass in the introduction we get
$$w\biggl(\pmatrix{*&*\cr c^2&a},\pmatrix{1&0\cr-c^2&1}\biggl)=0$$
and hence
$$\Bigl\{{c^2(1-a)\atop a}\Bigr\}=\Bigl\{{c^2\atop a}\Bigr\}.$$
So we obtain
$$\Bigl[{c\atop a}\Bigr]=\Bigl\{{c^2\atop a}\Bigr\}$$
and moreover
$$\Bigl[{c_1c_2\atop a}\Bigr]=\Bigl\{{c_1^2c_2^2\atop a}\Bigr\}=
\Bigl[{c_1\atop a}\Bigr]\Bigl\{{c_2^2\atop a}\Bigr\}
=\Bigl[{c_1\atop a}\Bigr]\Bigl[{c_2\atop a}\Bigr].$$
This finishes the proof of Proposition \MenS.\qed
\smallskip
The main result about Mennicke symbols is that they are trivial [BMS], Theorem 3.6. Hence we obtain now
the important result.
\proclaim
{Proposition}
{The multiplier system $v$ is identically one on all
$$\pmatrix{M&0\cr0&M'^{-1}}\quad\hbox{for}\ M\in\Gamma_1[q].$$}
Triio%
\finishproclaim
From Lemma \emaA\ follows now
$$\Bigl\{{c^2\atop d}\Bigr\}=1\quad\hbox{and}\quad
\Bigl\{{c_1\atop d}\Bigr\}=\Bigl\{{c_1c_2^2\atop d}\Bigr\}$$
for $c\equiv c_1\equiv c_2\equiv0$ mod $q$ and $d\equiv1$ mod $q$.
This can be generalized. We have to consider the Kronecker symbol $\bigl({c\over d}\bigr)$.
For its definition and properties we refer to [Di]. We will need it only for $c\ne 0$ and for
odd $d$.
We collect some properties
(always assuming this condition)
$$\Bigl({c_1c_2\over d}\Bigr)=
\Bigl({c_1\over d}\Bigr)\Bigl({c_2\over d}\Bigr),\quad
\Bigl({c\over d_1d_2}\Bigr)=
\Bigl({c\over d_1}\Bigr)\Bigl({c\over d_2}\Bigr).$$
Assume $d>0$ or $c_1c_2>0$. Then
$$\Bigl({c_1\over d}\Bigr)=\Bigl({c_2\over d}\Bigr)\quad\hbox{if}\quad c_1\equiv c_2\mod d.$$
Also the relation
$$\Bigl({c\over d_1}\Bigr)=\Bigl({c\over d_2}\Bigr)\quad\hbox{if}\quad
\cases{d_1\equiv d_2\mod c\ \hbox{and}\  c\equiv 0\mod4,\cr
d_1\equiv d_2\mod 4c\ \hbox{and}\ c\equiv 2\mod 4\cr}$$
is valid. Finally we mention
$$\Bigl({c\over -1}\Bigr)=\cases{1&for $c>0$,\cr -1&for $c<0$.\cr}$$
Since one of the rules demands $c\equiv0\mod 4$, we will from now on assume that
$q\equiv0\mod 4$.
\proclaim
{Proposition}
{Let $q$ be a suitable multiple of $q'$ and let
$$M=\pmatrix{a&b\cr c&d}\in\Gamma_1[q],\quad\Bigl({c\over d}\Bigr)=1.$$
Then $v(M)=1$.}
KronS%
\finishproclaim
{\it Proof.} We use the invariance under $(c,d)\mapsto (c,d+xc)$. We can apply Dirichlet's prime number
theorem and therefore assume that $d=p$ is a (positive) prime.
But then the Kronecker symbol is the usual
Legendre symbol.
Since $d\equiv 1\mod q$ we have $\bigl({q\over d}\bigr)=1$.
This implies $\bigl({c/q\over d}\bigr)=1$. Since $d$ is a prime, we get a solution
of $c/q=x^2+dy$ or $c=qx^2+dqy$.
Now use
$$\pmatrix{*&*\cr qx^2&d}\pmatrix{1&0\cr qy&1}=\pmatrix{*&*\cr c&d}.$$
In the case $c>0$ the $w$-value is zero. This follows from the corollary in the table of
Maass in the introduction.
In the case $c<0$ we must have $y<0$ and again from this corollary follows that the $w$-value is zero.
(In the notation of the table the sign distribution of $(m_1,c,m_1')$ is $(+,*,+)$ or $(+,-,*)$.)
Now we get
$$v(M)=v\pmatrix{*&*\cr c&d}=v\pmatrix{*&*\cr qx^2&d}=\Bigl\{{qx^2\atop d}\Bigr\}.$$
Lemma \emaA\ now shows
$$\Bigl\{{x^2q\atop d}\Bigr\}=\Bigl\{{x^2q^3\atop d}\Bigr\}=
\Bigl\{{q(qx)^2\atop d}\Bigr\}=\Bigl\{{q\atop d}\Bigr\}=\Bigl\{{q\atop 1}\Bigr\}=1.\eqno\square$$
\proclaim
{Lemma}
{Assume that the matrix $M={a\,b\choose c\,d}$ is contained in $\Gamma_1[q]$ and has
the following properties.
All entries are positive and $dq<c(q-1)$. Then
$$v(M)=e^{-2\pii r}\quad\hbox{if}\quad \Bigl({c\over d}\Bigr)=-1.$$}
zPir%
\finishproclaim
{\it Proof.} We consider
$$\pmatrix{a&b\cr c&d}\pmatrix{1-q&-q\cr q&1+q}=\pmatrix{*&*\cr c-qc+dq&-cq+d+dq}.$$
Clearly $\bigl({q\over 1+q}\bigr)=1$. We also claim
$$\Bigl({c-qc+dq\over -cq+d+dq}\Bigr)=1.$$
To prove this, we observe
$$\Bigl({c-qc+dq\over -cq+d+dq}\Bigr)=\Bigl({c-qc+dq\over d-c}\Bigr)=
\Bigl({c-qc+dq\over -1}\Bigr)\Bigl({c-qc+dq\over c-d}\Bigr).$$
Now we use $c-qc+dq<0$. It follows $c-d>0$. Hence we get
$$=-\Bigl({c\over c-d}\Bigr)=-\Bigl({c\over -d}\Bigr)=-\Bigl({c\over d}\Bigr)=-(-1)=1.$$
Now we have proved
$$v\biggl(\pmatrix{a&b\cr c&d}\pmatrix{1-q&-q\cr q&1+q}\biggr)=1.$$
The left hand side equals
$$v\pmatrix{a&b\cr c&d}\,\exp\left\{ 2\pii r w\left(\pmatrix{a&b\cr c&d},
\pmatrix{1-q&-q\cr q&1+q}\right)\right\}=1.$$
From Maass'  table in the introduction follows that the $w$-value is 1.
(The sign distribution of $(m_1,c,m_1')$ is $(+,+,-)$.)
This proves Lemma \zPir.\qed
\smallskip
There exist two coprime natural numbers
$c,d$ such that $c\equiv 0\mod q$ and $d\equiv 1\mod q$ and such that
$\bigl({c\over d}\bigl)=-1$. We also can assume $dq<c(q-1)$. The pair $(c,d)$ is the second row of
a matrix $M={a\,b\choose c\,d}\in\Gamma_1[q]$. We want to compute $v(M)$. Since we can add a multiple of the
second row to the first one, we can assume that $a$ and $b$ are also positive.
From Lemma \zPir\ we know $v(M)=e^{-2\pii r}$.
Now we consider
$$v(M^2)=v(M)^2e^{2\pii r w(M,M)}.$$
Since all entries from $M$ are positive, we have $w(M,M)=0$. So we get
$$v(M^2)=e^{-4\pii r}.$$
We compute $\bigl({\gamma\over\delta}\bigl)$ for the matrix
$$N=M^2=\pmatrix{\alpha&\beta\cr\gamma&\delta}.$$
We get
$$\Bigl({\gamma\over\delta}\Bigr)=\Bigl({c(a+d)\over cb+d^2}\Bigr)=
\Bigl({c\over cb+d^2}\Bigr)\Bigl({a+d\over cb+d^2}\Bigr).$$
We have
$$\Bigl({c\over cb+d^2}\Bigr)=\Bigl({c\over d^2}\Bigl)=1$$
and
$$\Bigl({a+d\over cb+d^2}\Bigr)=\Bigl({a+d\over d(a+d)-1}\Bigr).$$
Since $a+d\equiv 2\mod 4$ we only can change the denominator mod $4(a+d)$. Since $d\equiv1$ mod 4 we see
$$\Bigl({a+d\over d(a+d)-1}\Bigr)=\Bigl({a+d\over a+d-1}\Bigr)=\Bigl({1\over a+d-1}\Bigr)=1.$$
This shows $v(N)=1$ and we get  the relation
$$e^{-4\pii r}=1$$
which implies that $2r$ is integral. This finishes the proof of the main result.
\vfill\eject\noindent
{\paragratit References}
\bigskip
\item{[BMS]} Bass, H. Milnor, J. Serre, J.P.: {\it Solution of the congruence subgroup problem
for $\SL_n\; (n\ge 3)$ and $\Sp_{2n}\; (n\ge 2)$,}
Publications math\`ematiques l'I.H.\'E.S., tome {\bf 33}, p. 59--137 (1967)
\medskip
\item{[Ch]} Christian, U.:
{\it Hilbert-Siegelsche Modulformen und Poincar\'esche Reihen,} Mathematische
Annalen {\bf148}, 257–307  (1962)
\medskip
\item{[De]} Deligne, P.: {\it Extensions centrales non r\'esiduellement finies de groupes
arithmetiques,} C. R. Acad. Sci. Paris {\bf 287}, p. 203-208 (1978)
\medskip
\item{[Di]} Dickson, L.E.: {\it Introduction to the Theory of Numbers,}
Dover Publications, New York, Dover (1957).
\medskip
\item{[Hi]} Hill, R.: {\it Fractional weights and non-congruence subgroups,} Automorphic
Forms and Representations of algebraic groups over local fields,
Saito,\ H., Takahashi,\ T. (ed.) Surikenkoukyuroku series {\bf 1338}, 71-80 (2003)
\medskip
\item{[Ma1]} Maass, H.: {\it Lectures on Modular Functions of One Complex Variable,} Notes by Sunder
Lal, Tata Institute Of Fundamental Research, Bombay, Revised 1983  (1964)
\medskip
\item{[Ma2]} Maass, H.: {\it Die Multiplikatorsysteme zur Siegelschen Modulgruppe,}
Nach\-rich\-ten der Akademie der Wissenschaften zu G"ottingen
II. Mathematisch-physikalische Klasse, Nr. {\bf11}, 125-135 (1964)
\medskip
\item{[Me]} Mennicke, J.: {\it Zur Theorie der Siegelschen Modulgruppe,}
Math.\ Annalen {\bf 159}, 115--129 (1965)
\bye

%% file: neumac.tex
\input german

%
\output={\if N\header\headline={\hfill}\fi
\plainoutput\global\let\header=Y}
\magnification\magstep1
\tolerance = 500
\hsize=14.4true cm
\vsize=22.5true cm
\parindent=6true mm\overfullrule=2pt
\newcount\kapnum \kapnum=0
\newcount\parnum \parnum=0
\newcount\procnum \procnum=0
\newcount\nicknum \nicknum=1
\font\ninett=cmtt9

\font\ninebf=cmbx9

\font\sixbf=cmbx6
\font\ninesl=cmsl9

\font\nineit=cmti9

\font\ninerm=cmr9

\font\sixrm=cmr6
\font\ninei=cmmi9
\font\eighti=cmmi8
\font\sixi=cmmi6
\skewchar\ninei='177 \skewchar\eighti='177 \skewchar\sixi='177
\font\ninesy=cmsy9
\font\eightsy=cmsy8
\font\sixsy=cmsy6
\skewchar\ninesy='60 \skewchar\eightsy='60 \skewchar\sixsy='60
\font\titelfont=cmr10 scaled 1440
\font\paragratit=cmbx10 scaled 1200

\font\name=cmcsc10
\font\emph=cmbxti10

\font\tenmsbm=msbm10
\font\sevenmsbm=msbm7
%

%

%
\font\teneufm=eufm10
\font\seveneufm=eufm7
\font\fiveeufm=eufm5
\newfam\eufmfam
\textfont\eufmfam=\teneufm
\scriptfont\eufmfam=\seveneufm
\scriptscriptfont\eufmfam=\fiveeufm

\font\tenmsam=msam10
\font\sevenmsam=msam7
\font\fivemsam=msam5
\newfam\msamfam
\textfont\msamfam=\tenmsam
\scriptfont\msamfam=\sevenmsam
\scriptscriptfont\msamfam=\fivemsam
\font\tenmsbm=msbm10
\font\sevenmsbm=msbm7
\font\fivemsbm=msbm5
\newfam\msbmfam
\textfont\msbmfam=\tenmsbm
\scriptfont\msbmfam=\sevenmsbm
\scriptscriptfont\msbmfam=\fivemsbm
\def\Bbb#1{{\fam\msbmfam\relax#1}}
\def\cz{{\kern0.4pt\Bbb C\kern0.7pt}
}
\def\ez{{\kern0.4pt\Bbb E\kern0.7pt}
}
\def\fz{{\kern0.4pt\Bbb F\kern0.3pt}}
\def\gz{{\kern0.4pt\Bbb Z\kern0.7pt}}
\def\hz{{\kern0.4pt\Bbb H\kern0.7pt}
}
\def\kz{{\kern0.4pt\Bbb K\kern0.7pt}
}
\def\nz{{\kern0.4pt\Bbb N\kern0.7pt}
}
\def\oz{{\kern0.4pt\Bbb O\kern0.7pt}
}
\def\rz{{\kern0.4pt\Bbb R\kern0.7pt}
}
\def\sz{{\kern0.4pt\Bbb S\kern0.7pt}
}
\def\pz{{\kern0.4pt\Bbb P\kern0.7pt}
}
\def\qz{{\kern0.4pt\Bbb Q\kern0.7pt}
}
\newskip\ttglue
\def\ninepoint{\def\rm{\fam0\ninerm}%
  \textfont0=\ninerm \scriptfont0=\sixrm \scriptscriptfont0=\fiverm
  \textfont1=\ninei \scriptfont1=\sixi \scriptscriptfont1=\fivei
  \textfont2=\ninesy \scriptfont2=\sixsy \scriptscriptfont2=\fivesy
  \textfont3=\tenex \scriptfont3=\tenex \scriptscriptfont3=\tenex
  \def\it{\fam\itfam\nineit}%
  \textfont\itfam=\nineit
  \def\sl{\fam\slfam\ninesl}%
  \textfont\slfam=\ninesl
  \def\bf{\fam\bffam\ninebf}%
  \textfont\bffam=\ninebf \scriptfont\bffam=\sixbf
   \scriptscriptfont\bffam=\fivebf
  \def\tt{\fam\ttfam\ninett}%
  \textfont\ttfam=\ninett
  \tt \ttglue=.5em plus.25em minus.15em
  \normalbaselineskip=11pt
  \font\name=cmcsc9
  \let\sc=\sevenrm
  \let\big=\ninebig
  \setbox\strutbox=\hbox{\vrule height8pt depth3pt width0pt}%
  \normalbaselines\rm
  \def\sl{\it}}

\headline={\ifodd\pageno\rightheadline\else\leftheadline\fi}
\def\rightheadline{\ninepoint Paragraphen"uberschrift\hfill\folio}
\def\leftheadline{\ninepoint\folio\hfill Chapter"uberschrift}
\let\header=Y
\def\titel#1{\need 9cm \vskip 2truecm
\parnum=0\global\advance \kapnum by 1
{\baselineskip=16pt\lineskip=16pt\rightskip0pt
plus4em\spaceskip.3333em\xspaceskip.5em\pretolerance=10000\noindent
\titelfont Chapter \uppercase\expandafter{\romannumeral\kapnum}.
#1\vskip2true cm}\def\leftheadline{\ninepoint
\folio\hfill Chapter \uppercase\expandafter{\romannumeral\kapnum}.
#1}\let\header=N
}
\def\Titel#1{\need 9cm \vskip 2truecm
\global\advance \kapnum by 1
{\baselineskip=16pt\lineskip=16pt\rightskip0pt
plus4em\spaceskip.3333em\xspaceskip.5em\pretolerance=10000\noindent
\titelfont\uppercase\expandafter{\romannumeral\kapnum}.
#1\vskip2true cm}\def\leftheadline{\ninepoint
\folio\hfill\uppercase\expandafter{\romannumeral\kapnum}.
#1}\let\header=N
}
\def\need#1cm {\par\dimen0=\pagetotal\ifdim\dimen0<\vsize
\global\advance\dimen0by#1 true cm
\ifdim\dimen0>\vsize\vfil\eject\noindent\fi\fi}
\def\neupara#1{\par\penalty-2000
\procnum=0\global\advance\parnum by 1
\vskip1cm\noindent{\paragratit \the\parnum. #1}%
\def\rightheadline{\ninepoint\S\the\parnum.\ #1\hfill \folio}%
\vskip 8mm\noindent}
\def\Proclaim #1 #2\finishproclaim {\bigbreak\noindent
{\bf#1\unskip{}. }{\it#2}\medbreak\noindent}
%
\gdef\proclaim #1 #2 #3\finishproclaim {\bigbreak\noindent%
\global\advance\procnum by 1
{%
{\relax\ifodd \nicknum
\hbox to 0pt{\vrule depth 0pt height0pt width\hsize
   \quad \ninett#3\hss}\else {}\fi}%
\bf\the\parnum.\the\procnum\ #1\unskip{}. }
{\it#2}
\immediate\write\num{\string\def
 \expandafter\string\csname#3\endcsname
 {\the\parnum.\the\procnum}}
\medbreak\noindent}
\newcount\stunde \newcount\minute \newcount\hilfsvar
\def\uhrzeit{
    \stunde=\the\time \divide \stunde by 60
    \minute=\the\time
    \hilfsvar=\stunde \multiply \hilfsvar by 60
    \advance \minute by -\hilfsvar
    \ifnum\the\stunde<10
    \ifnum\the\minute<10
    0\the\stunde:0\the\minute~Uhr
    \else
    0\the\stunde:\the\minute~Uhr
    \fi
    \else
    \ifnum\the\minute<10
    \the\stunde:0\the\minute~Uhr
    \else
    \the\stunde:\the\minute~Uhr
    \fi
    \fi
    }

\def\Arg{\mathop{\rm Arg}\nolimits}

\def\GL{\mathop{\rm GL}\nolimits}

\def\im{\mathop{\rm Im}\nolimits}

\def\mod{\mathop{\rm mod}\nolimits}

\def\SL{\mathop{\rm SL}\nolimits}

\def\Sp{\mathop{\rm Sp}\nolimits}

\def\boxit#1{
  \vbox{\hrule\hbox{\vrule\kern6pt
  \vbox{\kern8pt#1\kern8pt}\kern6pt\vrule}\hrule}}
\def\Boxit#1{
  \vbox{\hrule\hbox{\vrule\kern2pt
  \vbox{\kern2pt#1\kern2pt}\kern2pt\vrule}\hrule}}

\def\smallni{\smallskip\noindent }

\def\lo{\longrightarrow}

\def\imag{{\rm i}}
\def\pii{\pi {\rm i}}

\def\square{\hbox{\hbox to 0pt{$\sqcup$\hss}\hbox{$\sqcap$}}}
\def\qed{\ifmmode\square\else{\unskip\nobreak\hfil
\penalty50\hskip3em\null\nobreak\hfil\square
\parfillskip=0pt\finalhyphendemerits=0\endgraf}\fi}
\def\pn{\the\parnum.\the\procnum}
\def\downmapsto{{\buildrel
        {\vbox{\hbox{\hskip.2pt$\scriptstyle-$}}}
        \over{\raise7pt\vbox{\vskip-4pt\hbox{$\textstyle\downarrow$}}}}}

%% file: german.tex
\expandafter\ifx\csname mdqon\endcsname\relax
\else   \fi

\message{Document Style Option `german'  Version 2 as of 16 May 1988}

\ifx\protect\undefined
\let\protect=\relax \fi

{\catcode`\@=11 

\gdef\allowhyphens{\penalty\@M \hskip\z@skip}

\newcount\U@C\newbox\U@B\newdimen\U@D
\gdef\umlauthigh{\def\"##1{{\accent127 ##1}}}
\gdef\umlautlow{\def\"{\protect\newumlaut}}
\gdef\newumlaut#1{\leavevmode\allowhyphens
     \vbox{\baselineskip\z@skip \lineskip.25ex
     \ialign{##\crcr\hidewidth
     \setbox\U@B\hbox{#1}\U@D .01\p@\U@C\U@D
     \U@D\ht\U@B\advance\U@D -1ex\divide\U@D \U@C
     \U@C\U@D\U@D\the\fontdimen1\the\font
     \multiply\U@D \U@C\divide\U@D 100\kern\U@D
     \vbox to .20ex  
     {\hbox{\char127}\vss}\hidewidth\crcr#1\crcr}}\allowhyphens}
\gdef\highumlaut#1{\leavevmode\allowhyphens
     \accent127 #1\allowhyphens}

\gdef\set@low@box#1{\setbox1\hbox{,}\setbox\z@\hbox{#1}\dimen\z@\ht\z@
     \advance\dimen\z@ -\ht1
     \setbox\z@\hbox{\lower\dimen\z@ \box\z@}\ht\z@\ht1 \dp\z@\dp1}

\gdef\@glqq{{\ifhmode \edef\@SF{\spacefactor\the\spacefactor}\else
     \let\@SF\empty \fi \leavevmode
     \set@low@box{''}\box\z@\kern-.04em\allowhyphens\@SF\relax}}
\gdef\glqq{\protect\@glqq}
\gdef\@grqq{\ifhmode \edef\@SF{\spacefactor\the\spacefactor}\else
     \let\@SF\empty \fi \kern-.07em``\kern.07em\@SF\relax}
\gdef\grqq{\protect\@grqq}
\gdef\@glq{{\ifhmode \edef\@SF{\spacefactor\the\spacefactor}\else
     \let\@SF\empty \fi \leavevmode
     \set@low@box{'}\box\z@\kern-.04em\allowhyphens\@SF\relax}}
\gdef\glq{\protect\@glq}
\gdef\@grq{\kern-.07em`\kern.07em}
\gdef\grq{\protect\@grq}
\gdef\@flqq{\ifhmode \edef\@SF{\spacefactor\the\spacefactor}\else
     \let\@SF\empty \fi
     \ifmmode \ll \else \leavevmode
     \raise .2ex \hbox{$\scriptscriptstyle \ll $}\fi \@SF\relax}
\gdef\flqq{\protect\@flqq}
\gdef\@frqq{\ifhmode \edef\@SF{\spacefactor\the\spacefactor}\else
     \let\@SF\empty \fi
     \ifmmode \gg \else \leavevmode
     \raise .2ex \hbox{$\scriptscriptstyle \gg $}\fi \@SF\relax}
\gdef\frqq{\protect\@frqq}
\gdef\@flq{\ifhmode \edef\@SF{\spacefactor\the\spacefactor}\else
     \let\@SF\empty \fi
     \ifmmode < \else \leavevmode
     \raise .2ex \hbox{$\scriptscriptstyle < $}\fi \@SF\relax}
\gdef\flq{\protect\@flq}
\gdef\@frq{\ifhmode \edef\@SF{\spacefactor\the\spacefactor}\else
     \let\@SF\empty \fi
     \ifmmode > \else \leavevmode
     \raise .2ex \hbox{$\scriptscriptstyle > $}\fi \@SF\relax}
\gdef\frq{\protect\@frq}

\global\let\original@ss=\ss
\gdef\newss{\leavevmode\allowhyphens\original@ss\allowhyphens}

\global\let\ss=\newss
\global\let\original@three=\3 


\gdef\german@dospecials{\do\ \do\\\do\{\do\}\do\$\do\&%
  \do\#\do\^\do\^^K\do\_\do\^^A\do\%\do\~\do\"}

\gdef\german@sanitize{\@makeother\ \@makeother\\\@makeother\$\@makeother\&%
\@makeother\#\@makeother\^\@makeother\^^K\@makeother\_\@makeother\^^A%
\@makeother\%\@makeother\~\@makeother\"}

\global\let\original@dospecials\dospecials
\global\let\dospecials\german@dospecials
\global\let\original@sanitize\@sanitize
\global\let\@sanitize\german@sanitize

\gdef\mdqon{\let\dospecials\german@dospecials
        \let\@sanitize\german@sanitize\catcode`\"\active}
\gdef\mdqoff{\catcode`\"12\let\original@dospecials\dospecials
        \let\@sanitize\original@sanitize}

{\mdqoff
\gdef\@UMLAUT{\"}
\gdef\@MATHUMLAUT{\mathaccent"707F}
\gdef\@SS{\mathchar"7019}
\gdef\dq{"}
}

{\mdqon
\gdef"#1{\if\string#1`\glqq{}%
\else\if\string#1'\grqq{}%
\else\if\string#1a\ifmmode\@MATHUMLAUT a\else\@UMLAUT a\fi
\else\if\string#1o\ifmmode\@MATHUMLAUT o\else\@UMLAUT o\fi
\else\if\string#1u\ifmmode\@MATHUMLAUT u\else\@UMLAUT u\fi
\else\if\string#1A\ifmmode\@MATHUMLAUT A\else\@UMLAUT A\fi
\else\if\string#1O\ifmmode\@MATHUMLAUT O\else\@UMLAUT O\fi
\else\if\string#1U\ifmmode\@MATHUMLAUT U\else\@UMLAUT U\fi
\else\if\string#1e\ifmmode\@MATHUMLAUT e\else\protect \highumlaut e\fi
\else\if\string#1i\ifmmode\@MATHUMLAUT i\else\protect\highumlaut\i \fi
\else\if\string#1E\ifmmode\@MATHUMLAUT E\else\protect\highumlaut E\fi
\else\if\string#1I\ifmmode\@MATHUMLAUT I\else\protect\highumlaut I\fi
\else\if\string#1s\ifmmode\@SS\else\ss\fi{}%
\else\if\string#1-\allowhyphens\-\allowhyphens
\else\if\string#1\string"\hskip\z@skip
\else\if\string#1|\discretionary{-}{}{\kern.03em}%
\else\if\string#1c\allowhyphens\discretionary{k-}{}{c}\allowhyphens
\else\if\string#1f\allowhyphens\discretionary{ff-}{}{f}\allowhyphens
\else\if\string#1l\allowhyphens\discretionary{ll-}{}{l}\allowhyphens
\else\if\string#1m\allowhyphens\discretionary{mm-}{}{m}\allowhyphens
\else\if\string#1n\allowhyphens\discretionary{nn-}{}{n}\allowhyphens
\else\if\string#1p\allowhyphens\discretionary{pp-}{}{p}\allowhyphens
\else\if\string#1t\allowhyphens\discretionary{tt-}{}{t}\allowhyphens
\else\if\string#1<\flqq{}%
\else\if\string#1>\frqq{}%
\else             \dq #1%
\fi\fi\fi\fi\fi\fi\fi\fi\fi\fi\fi\fi\fi\fi\fi\fi\fi\fi\fi\fi\fi\fi\fi\fi\fi}
}

\gdef\dateaustrian{\def\today{\number\day.~\ifcase\month\or
  J\"anner\or Februar\or M\"arz\or April\or Mai\or Juni\or
  Juli\or August\or September\or Oktober\or November\or Dezember\fi
  \space\number\year}}
\gdef\dategerman{\def\today{\number\day.~\ifcase\month\or
  Januar\or Februar\or M\"arz\or April\or Mai\or Juni\or
  Juli\or August\or September\or Oktober\or November\or Dezember\fi
  \space\number\year}}
\gdef\dateUSenglish{\def\today{\ifcase\month\or
 January\or February\or March\or April\or May\or June\or
 July\or August\or September\or October\or November\or December\fi
 \space\number\day, \number\year}}
\gdef\dateenglish{\def\today{\ifcase\day\or
 1st\or 2nd\or 3rd\or 4th\or 5th\or
 6th\or 7th\or 8th\or 9th\or 10th\or
 11th\or 12th\or 13th\or 14th\or 15th\or
 16th\or 17th\or 18th\or 19th\or 20th\or
 21st\or 22nd\or 23rd\or 24th\or 25th\or
 26th\or 27th\or 28th\or 29th\or 30th\or
 31st\fi
 ~\ifcase\month\or
 January\or February\or March\or April\or May\or June\or
 July\or August\or September\or October\or November\or December\fi
 \space \number\year}}
\gdef\datefrench{\def\today{\ifnum\day=1\relax 1\/$^{\rm er}$\else
  \number\day\fi \space\ifcase\month\or
  janvier\or f\'evrier\or mars\or avril\or mai\or juin\or
  juillet\or ao\^ut\or septembre\or octobre\or novembre\or d\'ecembre\fi
  \space\number\year}}


\gdef\captionsgerman{%
\def\refname{Literatur}%
\def\abstractname{Zusammenfassung}%
\def\bibname{Literaturverzeichnis}%
\def\chaptername{Kapitel}%
\def\appendixname{Anhang}%
\def\contentsname{Inhaltsverzeichnis}%
\def\listfigurename{Abbildungsverzeichnis}%
\def\listtablename{Tabellenverzeichnis}%
\def\indexname{Index}%
\def\figurename{Abbildung}%
\def\tablename{Tabelle}%
\def\partname{Teil}}

\gdef\captionsenglish{%
\def\refname{References}%
\def\abstractname{Abstract}%
\def\bibname{Bibliography}%
\def\chaptername{Chapter}%
\def\appendixname{Appendix}%
\def\contentsname{Contents}%
\def\listfigurename{List of Figures}%
\def\listtablename{List of Tables}%
\def\indexname{Index}%
\def\figurename{Figure}%
\def\tablename{Table}%
\def\partname{Part}}

\gdef\captionsfrench{%
\def\refname{R\'ef\'erences}%
\def\abstractname{R\'esum\'e}%
\def\bibname{Bibliographie}%
\def\chaptername{Chapitre}%
\def\appendixname{Appendice}%
\def\contentsname{Table des mati\`eres}%
\def\listfigurename{Liste des figures}%
\def\listtablename{Liste des tables}%
\def\indexname{Index}%
\def\figurename{Figure}%
\def\tablename{Table}%
\def\partname{Partie}}%

\newcount\language 
\newcount\USenglish  \global\USenglish=0
\newcount\german     \global\german=1
\newcount\austrian   \global\austrian=2
\newcount\french     \global\french=3
\newcount\english    \global\english=4
\gdef\setlanguage#1{\language #1\relax
  \expandafter\ifcase #1\relax
  \dateUSenglish  \captionsenglish   \or
  \dategerman     \captionsgerman    \or
  \dateaustrian   \captionsgerman    \or
  \datefrench     \captionsfrench    \or
  \dateenglish    \captionsenglish   \fi}

\gdef\originalTeX{\mdqoff \umlauthigh
  \let\ss\original@ss \let\3\original@three
  \setlanguage{\USenglish}}
\gdef\germanTeX{\mdqon \umlautlow \let\ss\newss \let\3\ss
  \let\dospecials\german@dospecials
  \setlanguage{\german}}

} 


\germanTeX
